\documentclass{article}
\usepackage[english]{babel}
\usepackage[margin=1in]{geometry}
\usepackage{multicol}
\usepackage[math]{blindtext}

\usepackage{url}
\usepackage{amssymb,amsmath}
\usepackage{amsmath}
\usepackage{amssymb}
\usepackage{dsfont}
\usepackage{varioref,psfrag}
\usepackage{multirow}
\usepackage{indentfirst}       
\usepackage{graphicx} 
\usepackage{subfig}  
\usepackage{url}              
\usepackage{textcomp}
\usepackage{xspace}
\usepackage{cite}
\usepackage{filecontents}
\usepackage{algorithm}
\usepackage{algpseudocode}
\usepackage{color}
\usepackage{footnote}
\usepackage{float}
\restylefloat{table}
\usepackage{color}

\newcommand{\RNum}[1]{\uppercase\expandafter{\romannumeral #1\relax}}

\begin{document}

{\centering

{\bfseries\Large Decentralized Clustering and Linking  by Networked
Agents \bigskip}

Sahar Khawatmi\textsuperscript{1}, Ali H.   Sayed\textsuperscript{2}, and
 Abdelhak M. Zoubir\textsuperscript{1} \\
  {\itshape
\textsuperscript{1}	Technische Universit\"at
 Darmstadt, Signal Processing Group, 
		  Darmstadt, Germany \\ 
\textsuperscript{2} University of
  California,  Department of Electrical  Engineering,
		 Los Angeles, USA\\
  }
  \vspace{0.5cm}
\centerline{October 7, 2016}}

\begin{abstract}
We consider the problem of decentralized clustering and estimation over
multi-task networks, where agents infer and track
different models of interest.
The agents do not know beforehand  which model is generating their own data.
 They also do not know which agents in their
neighborhood belong to the same cluster.
We propose a decentralized clustering algorithm aimed 
at identifying and forming clusters of 
agents of similar objectives, and at guiding cooperation
 to enhance the inference performance. 
 One key feature of the proposed technique is the integration
  of the learning and clustering tasks into a single strategy.
   We analyze the performance of the procedure and show
    that the error probabilities of types I and II decay exponentially 
    to zero with the step-size parameter.
While links between agents following different objectives 
are ignored in the clustering process, 
we nevertheless show how to exploit these 
links to relay critical information across the network 
for enhanced performance. Simulation results illustrate the 
performance of the proposed method in comparison to other useful techniques.
  \let\thefootnote\relax\footnote{The work of Ali. H. Sayed was supported
  in part by NSF grants ECCS--1407712 and CCF--1524250.
  An early short version of this work appeared in the conference
  publication~\cite{cluster_me}. This article has been submitted for publication.} 
\end{abstract}

\section{Introduction and Related Work}

\noindent Distributed learning is a powerful technique for extracting
information from networked agents 
(see, e.g.,  ~\cite{adapt2,book2,Poor,adapt3,ref13,ref8,al2016robust} and
the references therein).  In this work, we consider  a network of agents connected
by a graph.
Each agent senses data
generated by some {\em unknown} model.
It is assumed that there are clusters of agents within the network,
where agents in the same cluster observe data arising from the
same model. 

However, the   agents do not know which model is generating their own data.
They also do not know which agents in their neighborhood  
belong to the same cluster. Scenarios of this type arise, for example,
in tracking applications when a collection of networked 
agents is tasked with tracking several moving
objects~\cite{ref1,ref2,mobile1}.
Clusters end up being formed within the network with 
different clusters following different targets. 
The quality of the tracking/estimation performance
will be improved if neighboring agents following the 
same target know of each other to promote cooperation. 
It is not only cooperation within clusters that is useful,
but also cooperation across clusters, especially 
when targets move in formation and the location of the 
targets  are correlated. 
Motivated by these considerations,
the main objective of this work is to develop 
a distributed technique that enables agents to recognize 
neighbors from the same cluster and promotes cooperation 
for improved inference performance.

There have been several  useful works in the 
literature on the solution of inference problems for such  multi-task networks, i.e., for
networks with multiple unknown models (also called tasks) --- see,
e.g.,~\cite{ref3,ref6,dm,multi1,adapt3,adapt4,cluster1,cluster2,cluster4} and
the references therein.
In the solutions developed in~\cite{cluster1,cluster2,cluster4}, 
clustering is achieved by relying on adaptive combination strategies, 
whereby   weights on  edges between agents 
are adapted and their size becomes smaller for unrelated tasks. 
In  these earlier works, there still exists the possibility that valid 
links between agents belonging to the same cluster may be overlooked, 
mainly due to errors during the adaptation process. 
A more robust clustering method was proposed in~\cite{cluster3} where the
clustering and learning operations were decoupled from each other. 
In this way, tracking errors do not influence the clustering 
mechanism and the resulting distributed algorithm enables 
the agents to identify their clusters and to attain improved learning 
 accuracy. The work~\cite{cluster3} evaluated the error probabilities of types I
 and II, i.e., of false  alarm and mis-detection 
 for their proposed scheme and showed that these errors
  decay exponentially 
   with the step-size.
   This means that, under their  scheme, the probability of correct
   clustering can be made arbitrarily close to one by selecting sufficiently small step-sizes.

Still, it is preferable to {\em merge}  the clustering and 
learning mechanisms rather than have them run separately of each other. 
Doing so reduces the computational burden and, if successful, 
can also lead to  enhancement in clustering accuracy 
relative to the earlier approaches~\cite{cluster1,cluster2,cluster4}.
 We showed
in preliminary work~\cite{cluster_me} that this is indeed possible 
for a particular class of inference problems involving  mean-square-error risks. 
In this work, we  generalize the results and devise an {\em integrated}  
clustering-learning approach for general-purpose risk functions. 
Additionally, and  motivated by the results from~\cite{dm}
on adaptive decision-making by networked agents, 
we further  incorporate a smoothing mechanism  into our strategy to enhance 
  the belief that agents have about their clusters. 
We also  show how to exploit the unused links 
among neighboring agents belonging to different clusters to 
{\em relay}  useful information among agents. 
We carry out a detailed analysis of the resulting framework, and illustrate its superior 
performance  by means of 
computer  simulations.

The organization of the work is as follows. The network and data model are
described in Section \RNum{2}, while the integrated clustering and inference
framework is developed  in Section \RNum{3}. The network error recursions are derived
in Section \RNum{4}, and the  probabilities of erroneous decision are derived  in Section \RNum{5}.
In Section \RNum{6} we illustrate the linking technique for relaying information, 
and present simulation results in Section \RNum{7}.\\
    
\noindent {\bf Notation}. We use lowercase letters to denote vectors, uppercase
letters for matrices, plain letters for deterministic variables, and boldface
letters for random variables. The superscript $\circ$ is used to indicate true
values. The letter $\mathbb{E}$ denotes the expectation operator. The Euclidean norm is denoted by 
$\|\cdot\|$. The symbols $\mathds{1}$ and $I$ denote
the all-one vector and the  identity matrix of appropriate sizes, respectively. We
write $(\cdot)^\intercal$, $(\cdot)^{-1}$, and $\textrm{Tr}(\cdot)
$ to denote transposition, matrix
 inversion, and matrix trace, respectively. The  $\textrm{diag}\{\cdot\}$
 operator extracts the diagonal entries of its matrix argument and stacks them into a
 column.
 The $k-$th row
 (column) of matrix $X$ is denoted by $[X]_{k,:}$ ($[X]_{:,k}$).

\section{Network and Data Model}

\subsection{Network  Overview}
\noindent  We consider a network with $N$ agents connected by a graph. 
It is assumed that there are $C$ clusters, denoted by 
${\cal C}_1, {\cal C}_2, \ldots, {\cal C}_C$, where each 
${\cal C}_m$ represents the set of agent indices in that cluster. 
We associate an unknown column vector of size $M\times 1$ with each cluster, 
denoted by $w^\circ_{{\cal C}_m}\in\mathbb{R}^{M}$.
 The aggregation of all these unknowns is denoted by
  \begin{equation}\label{eq:wo}
w^\circ_{\mathcal{C}}\triangleq\textrm{col}
\{w^\circ_{\mathcal{C}_{1}},w^\circ_{\mathcal{C}_{2}},\ldots,{w}^\circ_{\mathcal{C}_{C}}\}
, \ (CM\times1).
\end{equation}
Each agent $k$ wishes to recover  the model for its cluster; 
the unknown model for agent $k$ is  denoted by $\{w_k^\circ\}$.
Obviously, this model agrees with 
the model of the cluster that $k$ belongs to, i.e.,
 $w_k^\circ=w_{{\cal C}_m}^\circ$ if $k\in {\cal C}_m$. We stack all
 $\{w_k^\circ\}$ into a column vector:
\begin{equation}\label{eq:unknown_models}
{w}^\circ\triangleq\textrm{col}
\{{w}^\circ_{1},{w}^\circ_{2},\ldots,{w}^\circ_{N}\}
, \ (NM\times1).
\end{equation}

\begin{figure}
\centering\includegraphics[width=17pc]{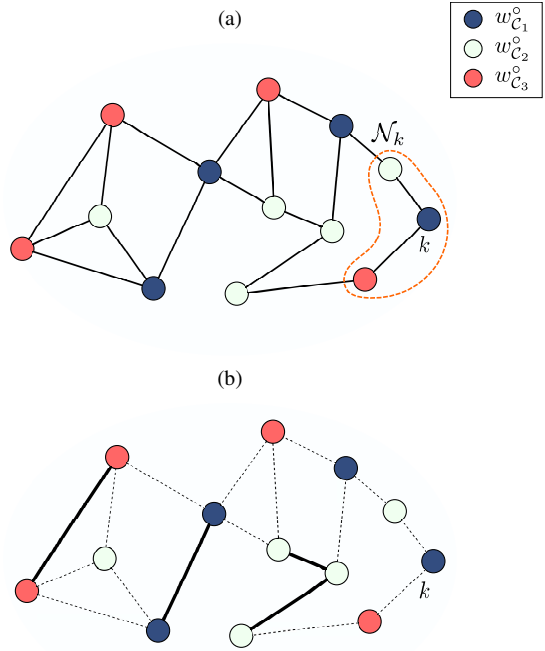}
\caption{\small ({\em Top}) Example  of a network topology involving
three clusters, represented by three different colors. 
({\em Bottom}) Clustered topology  that will result for the network shown on top.}
\label{fig:Partition0}
\end{figure}

Figure~\ref{fig:Partition0} illustrates a network with $C=3$ 
clusters represented by three colors. 
All agents in the same cluster are interested in
estimating the same parameter vector.  We denote the set of neighbors of an
agent $k$ by ${\cal N}_k$. Observe in this example that the neighbors of agent
$k$ belong to different clusters. The cluster information is 
not known to the agents beforehand. 
For instance, agent $k$ would not know that 
its neighbors are sensing data arising from different models. 
If we allow the network to perform indiscriminate cooperation, 
then performance will degrade significantly.  
For this reason, a clustering operation is needed to allow
the agents to learn which neighbors to cooperate with towards 
the same objective. The technique developed in this work will 
allow agents to emphasize links to neighbors 
in the same cluster and to disregard links to neighbors from other clusters.
 The outcome would be a graph structure similar to the one shown 
 in the bottom part of the same figure, where unwarranted links are shown in dotted lines.
 In this way, the interference caused by
different objectives is avoided and the overall performance for each
cluster will be  improved. Turning off a link  between two agents means that there is
no more sharing of data between them. Still, we will exploit 
these ``unused'' links by assigning to them a useful role 
in relaying information across the network.  

 \subsection{Topology Matrices}
\noindent In preparation for the description of the proposed strategy, 
we introduce the $N\times N$ adjacency matrix $Y=[y_{\ell k}]$, 
whose elements are either zero or one depending on whether
 agents are linked by an edge or not. Specifically, 
\begin{equation}
{y}_{\ell k}=\begin{cases} 
1, &  \ell \in \mathcal{N}_k, \\
0 , &  \textrm{otherwise}.
\end{cases}
\end{equation}
We assume that each  agent $k$ belongs to its neighborhood set, 
$k \in {\cal N}_k$. The set ${\cal N}^-_k$ excludes $k$.   
Agents know their neighborhoods but they do not know which subset 
of their neighbors is subjected to data from the same model.
 In order to devise a procedure that allows agents to
  arrive at this information, we introduce a second  $N\times N$ 
  clustering matrix, denoted by $E_i$ at time $i$, in a manner 
  similar to the  adjacency matrix  $Y$, except that the value at location
  $(\ell,k)$ will be set to one if agent $k$ {\em believes}  at time
   $i$ that its neighbor $\ell$ belongs to the same cluster:
   \begin{equation}
{e}_{\ell k}(i)=\begin{cases} 
1, &  \textrm{if } \ell\in{\cal N}_k \textrm{ and }  k  \textrm{ believes that
} w_k^\circ=w_\ell^\circ,\\
0 , &  \textrm{otherwise}.
\end{cases}
\end{equation}
   The entries of $E_i$ will be learned online. 
   At every time $i$, we can then use    these entries to infer which neighbors
   of $k$ are believed to belong to the same cluster as $k$;
    these would be the indices of the nonzero entries in the $k-$th column of
     $E_i$.
      We collect these indices into the neighborhood set,
       ${\cal N}_{k,i}$; this set is a subset of ${\cal N}_k$ and it 
       evolves over time during the learning process. 
       At any time $i$, agent $k$ will only be 
       cooperating with the neighbors within ${\cal N}_{k,i}$.
        We will describe in the sequel how $E_i$ is learned.
 
 \subsection{Problem Formulation}
\noindent We associate with each agent $k$  a strongly-convex  and differentiable risk
function $J_k(w_k)$, with a unique minimum  at
$w_k^\circ$.
In general, each risk $J_{k}(w_k): \mathbb{R}^{M}\rightarrow\mathbb{R}$, is defined as the
expectation of some loss function ${Q}_k(\cdot)$, say,
\begin{equation}
J_k(w_k)\;=\;\mathbb{E} \; {Q}_k(w_k; \boldsymbol{z}_k)
\end{equation}
where $\boldsymbol{z}_k$ denotes random data sensed by agent $k$ and the expectation
 is over the distribution of this data. The 
 network of agents is interested  in estimating the 
 minimizers of the following aggregate cost over the vectors $\{w_k\}$:
\begin{equation}\label{eq:cost}
   J^{\textrm{glob}}(w_1,w_2,\ldots, w_N)\triangleq\sum_{k=1}^N J_k(w_k).
\end{equation} 
Since agents  from different clusters do not share the same minimizers,  
the aggregate cost can be re-written as 
\begin{equation}\label{eq:cost2}
  J^{\textrm{glob}}(w_{{\cal C}_1},\ldots, w_{{\cal C}_C})
   \triangleq\sum_{m=1}^C \sum_{k\in {\cal C}_m} J_k(w_k)
\end{equation} 
where $w_k^\circ=w_{{\cal C}_m}^\circ$. 
We collect the gradient vectors of the risk functions 
across the network into the aggregate vector
\begin{equation}
\label{eq:gradient}
\nabla J({w})  \triangleq \textrm{col}\{\nabla
J_1({w}_1),\ldots, \nabla J_{N}({w}_{N})\}.
\end{equation}
These gradients will not be available
 in most cases since the distribution of the data 
 is not known to enable evaluation of the expectation of the loss functions. 
 In stochastic-gradient implementations, it is customary to replace the above 
 aggregate vector by the following approximation where the true gradients 
 of the risk functions are replaced by
\begin{equation}
\label{eq:hat_gradient}
 \widehat{\nabla J}({w}) \triangleq \textrm{col}\{
 \widehat{\nabla J_1}({w}_1),\ldots, 
  \widehat{\nabla J_{N}}({w}_{N})\}
\end{equation}
where each $\widehat{\nabla J}_k(w_k)$ is constructed from the gradient of the
respective loss function
\begin{equation}
\widehat{\nabla J}_k(w_k)\;=\;\nabla {Q}_k(w_k;\boldsymbol{z}_k)
\end{equation}
evaluated at the corresponding data point, $\boldsymbol{z}_k$.

\subsection{Assumptions}

\noindent   We list here the assumptions that are needed to drive the analysis.
 These assumptions are typical in the analysis of stochastic-gradient
 algorithms, and most of them are automatically satisfied by important cases of interest, 
 such as when the risk functions are quadratic or logistic --- see,
 e.g.,~\cite{book2,cluster3}. 
 
 We thus assume that each  individual cost function
 $J_k(w_k)$ is twice-differentiable and  $\tau_k-$strongly
 convex~\cite{book2,book5,book8}, for some $\tau_k>0$.
 We also require the gradient vector of $J_k(w_k)$ to be $\zeta_k-$Lipschitz, i.e.,
\begin{equation}\label{eq:cost_bound2}
\|\nabla J_k(w_{k_1})-\nabla J_k(w_{k_2})\|
\leq
\zeta_k \|w_{k_1}-w_{k_2}\|
\end{equation} 
for any $w_{k_1},w_{k_2}\in \mathbb{R}^{M}$. Thus,
the Hessian matrix function $\nabla^2 J_k(w_k)$  is bounded by
\begin{equation}\label{eq:cost_bounds}
\tau_k I_M\leq\nabla^2 J_k(w_k)\leq \ \zeta_k I_M
\end{equation} 
where $\tau_k \leq \zeta_k $. 
Each Hessian matrix function  is also assumed to satisfy the
Lipschitz condition:
\begin{equation}\label{eq:Lipschitz}
\|\nabla^2 J_k(w_{k_1})-\nabla^2 J_k(w_{k_2})\|
\leq
\kappa_k \|w_{k_1}-w_{k_2}\|
\end{equation}
for 
some $\kappa_k\geq0$ and any $w_{k_1},w_{k_2}\in \mathbb{R}^{M}$.
The network gradient noise is denoted by
$\boldsymbol{s}_i(\boldsymbol{w}_{i-1})$ which is the random process defined by 
\begin{equation}\label{eq:gradient_noise_vec}
\boldsymbol{s}_i(\boldsymbol{w}_{i-1}) \triangleq
\textrm{col}\{
 \boldsymbol{s}_{1,i}(\boldsymbol{w}_{1,i-1}),\ldots, 
 \boldsymbol{s}_{N,i}(\boldsymbol{w}_{N,i-1})\}
\end{equation}
where the  gradient noise at agent $k$ at  time  $i$  is given by
\begin{equation}\label{eq:gradient_noise}
\boldsymbol{s}_{k,i}(\boldsymbol{w}_{k,i-1})
 \triangleq\widehat{\nabla J}_k(\boldsymbol{w}_{k,i-1})- 
 {\nabla J}_k(\boldsymbol{w}_{k,i-1}).
\end{equation}
Here we are denoting the iterates $\boldsymbol{w}$ 
in boldface notation to indicate that
they will actually be stochastic variables
 due to the approximation of the true gradients.

We let  $\{\mathbb{F}_{k,i};i\geq0\}$ denote the filtration 
that collects all information up to time $i$. We then denote the conditional 
covariance matrix  of $\boldsymbol{s}_{k,i}(\boldsymbol{w}_{k,i-1})$ by
\begin{equation}\label{eq:gradient_noise_cov}
  { R}_{k,i}(\boldsymbol{w}_{k,i-1})\triangleq
   \mathbb{E}  \left[ \boldsymbol{s}_{k,i}(\boldsymbol{w}_{k,i-1})
   \boldsymbol{s}^{\intercal}_{k,i}(\boldsymbol{w}_{k,i-1}) \mid
   \mathbb{F}_{k,i-1}
   \right].
\end{equation}
It is assumed that the  gradient noise process satisfies the following properties for any 
$\boldsymbol{w}_{k,i-1}$ in $ \mathbb{F}_{k,i-1}$~\cite{book2}:
\begin{enumerate}
  \item \hspace{-0.3cm} Martingale difference~\cite{book2,cluster3}:
\begin{equation}\label{eq:gradient_noise_1}
   \mathbb{E}  \left[ \boldsymbol{s}_{k,i}(\boldsymbol{w}_{k,i-1})\mid
   \hspace{-0.1cm}  \mathbb{F}_{k,i-1}   \right]=0
\end{equation}
\item \hspace{-0.3cm} Bounded fourth-order moment~\cite{book2,cluster3}:
\begin{equation}\label{eq:gradient_noise_2}
   \mathbb{E}  \left[ \|\boldsymbol{s}_{k,i}(\boldsymbol{w}_{k,i-1})\|^4 \mid
  \hspace{-0.1cm} \mathbb{F}_{k,i-1}   \right]\leq
   \beta_k^2 \|w_k^\circ-\boldsymbol{w}_{k,i-1}\|^4+\rho_k^4
\end{equation}
for some $\beta^2_k,\rho^4_k\geq0$.
 \vspace{0.15cm}
\item \hspace{-0.3cm} Lipschitz conditional covariance
function~\cite{book2,cluster3}:
\begin{equation}\label{eq:gradient_noise_3}
 \|{{ R}}_{k,i}(w_k^\circ) -
 {{ R}}_{k,i}(\boldsymbol{w}_{k,i-1}) \|\leq
 \theta_k  \|w_k^\circ -\boldsymbol{w}_{k,i-1} \|^{\eta_k}
\end{equation}
for some $\theta_k\geq0$ and $0<\eta_k\leq4$.
 \vspace{0.15cm}
\item \hspace{-0.3cm} Convergent conditional covariance
matrix~\cite{book2,cluster3}:
\begin{equation}\label{eq:gradient_noise_4}
{ R}_k\triangleq
\lim_{i\rightarrow\infty} {{ R}}_{k,i}(w_k^\circ) > 0
\end{equation}
where ${R}_k$ is symmetric and positive definite.
 \end{enumerate}

 \subsection{Data Model}
 
 \noindent  We assume that each  agent $k$ runs an independent  stochastic gradient-descent algorithm of the form:
\begin{equation}\label{eq:data_step1}
 \boldsymbol{\psi}_{k,i}= \boldsymbol{\psi}_{k,i-1}- \mu_k 
 \widehat{\nabla J_k}(\boldsymbol{\psi}_{k,i-1})
\end{equation}
where $\mu_k>0$ is a small step-size parameter, and $\boldsymbol{\psi}_{k,i}$ denotes 
 the intermediate estimate for $w^\circ_k$ at time  $i$.
Cooperation among agents will be limited to neighbors that belong to the same
cluster. Therefore, following the update~(\ref{eq:data_step1}), ideally, agent
$k$ should 
only share data with agent $\ell$ if ${w}^\circ_{k}={w}^\circ_{\ell}$. 
The agents do not know which agents in their
neighborhood belong to the same cluster; this information is learned in real-time. Therefore, 
agent $k$ will only  share data with agent $\ell$ if it
believes that ${w}^\circ_{k}={w}^\circ_{\ell}$. Specifically, 
agent $k$ will combine the estimates 
from its neighbors in a convex manner as follows:
\begin{equation}\label{eq:data_step2}
  \boldsymbol{w}_{k,i}= \sum_{\ell=1}^{N}\boldsymbol{a}_{\ell k}(i)
 \boldsymbol{\psi}_{\ell,i}
\end{equation}
where the nonnegative combination coefficients $\{\boldsymbol{a}_{\ell k}(i)\}$ satisfy
\begin{equation}\label{eq:matrix_a}
\boldsymbol{a}_{kk}(i)>0, \   \boldsymbol{a}_{\ell k}(i) = 0 \ \ \textrm{for}  \ \ell \notin
\boldsymbol{\mathcal{N}}_{k,i}, \ \ \ \ \sum_{\ell=1}^N \boldsymbol{a}_{\ell
k}(i)=1.
\end{equation}
In the next section, we explain how the combination coefficients $\{\boldsymbol{a}_{\ell k}(i)\}$ are selected in
order to perform the combined tasks of estimation and clustering.

\section{Clustering Scheme}
\noindent Let $\delta>0$ denote the smallest distance among the cluster models,
$\{w^\circ_{{\cal
C}_m}\}$.  For  any distinct $a,b\in\{1,\ldots,C\}$, it then holds  that
 \begin{equation}\label{eq:delta_CR}
\|w^\circ_{\mathcal{C}_{a}}-w^\circ_{\mathcal{C}_{b}}\|\geq \delta.
 \end{equation}
  We  introduce an $N\times N$  trust matrix $\boldsymbol{F}_i$;
each entry $\boldsymbol{f}_{\ell k}(i)\in[0,1]$ of this matrix reflects the
amount of trust that agent $k$ has in neighbor $\ell \in {\cal
N}^-_{k}$ belonging to its cluster. 
The entries $\{\boldsymbol{f}_{\ell k}(i)\}$
are constructed as follows. Agent $k$ first computes the Boolean variable:
 \begin{equation}\label{eq:b}
 \boldsymbol{b}_{\ell k}(i)=\begin{cases}
 1, &   \textrm{if }\| \boldsymbol{\psi}_{\ell,i} - \boldsymbol{w}_{k,i-1}\|^2
 \leq \alpha,\\
 0,  & \textrm{otherwise}
 \end{cases}
 \end{equation}
  where $\alpha$ is a threshold value. The trust  level
  $\boldsymbol{f}_{\ell k}(i)$ is smoothed  as follows:
 \begin{equation}\label{eq:f}
   \boldsymbol{f}_{\ell k}(i)=\nu  \boldsymbol{f}_{\ell k}(i-1)+ (1-\nu) 
   \boldsymbol{b}_{\ell k}(i)
 \end{equation}
  where  the forgetting factor, $0<\nu<1$, determines the speed with which trust
  in neighbor $\ell$ accumulates over time.   Once $\boldsymbol{f}_{\ell k}(i)$
  exceeds some threshold  $\gamma$, agent $k$ declares that neighbor $\ell$ belongs to its 
   cluster and sets the corresponding entry $\boldsymbol{e}_{\ell k}(i)$ in
   matrix $\boldsymbol{E}_i$ to the value one:
  \begin{equation}\label{eq:e}
     \boldsymbol{e}_{\ell k}(i)=
  \begin{cases}
 1,\ \  \textrm{if }{ \boldsymbol{f}_{\ell k}(i) \geq \gamma},\\
 0, \ \   \textrm{otherwise}
  \end{cases}
 \end{equation}
 where $0<\gamma<1$. For completeness, we set  for any  agent $k$, 
 $\boldsymbol{b}_{kk}(i)=\boldsymbol{f}_{kk}(i)=\boldsymbol{e}_{kk}(i)=1$.
 Observe that the computation of the binary variable $\boldsymbol{b}_{\ell
 k}(i)$ couples the $\boldsymbol{\psi}$ and $\boldsymbol{w}$ variables.
 Therefore, by using smoothed values $\{\boldsymbol{f}_{\ell k}(i)\}$  for the trust variables, we are able to
  couple the clustering and inference procedures into a single iterative algorithm rather than run them separately. The smoothing  reduces the influence of erroneous clustering decisions on the inference task. 
    The following listing summarizes the proposed strategy.
 \begin{algorithm}[H]
  \begin{algorithmic}
  \State Initialize
  $\boldsymbol{F}_{-1}=\boldsymbol{B}_{-1}=\boldsymbol{E}_{-1}=I$ and
  $\boldsymbol{\psi}_{-1}=\boldsymbol{w}_{-1}=0$.
    \For {$i\geq 0$}
      \vspace{0.05cm}
   \For {$k=1,\ldots,N$}
    \vspace{-0.12cm}
       \begin{align}
   & \ \boldsymbol{\psi}_{k,i}= \boldsymbol{\psi}_{k,i-1}- \mu_k \widehat{\nabla
       J_k}(\boldsymbol{\psi}_{k,i-1})  \ &
  \end{align}
 \For {$\ell \in {\cal N}^-_{k}$}
  \vspace{-0.07cm}
 \begin{align}
 & \boldsymbol{b}_{\ell k}(i)=\begin{cases}
 1, &   \textrm{if } \| \boldsymbol{\psi}_{\ell,i} - \boldsymbol{w}_{k,i-1}\|^2
 \leq \alpha\\
 0,  & \textrm{otherwise}
 \end{cases}\\
 & \boldsymbol{f}_{\ell k}(i)=\nu  \boldsymbol{f}_{\ell k}(i-1)+ (1-\nu) 
 \; \boldsymbol{b}_{\ell k}(i)
\end{align}
\State \hspace{-0.35cm} update $\boldsymbol{e}_{\ell k}(i)$ according
to~(\ref{eq:e})
 \vspace{0.2cm}
 \EndFor
  \vspace{-0.07cm}
  \begin{align}
 &\textrm{select } \{\boldsymbol{a}_{\ell k}(i)\} \textrm{ according to}~(\ref{eq:matrix_a})
\textrm{ and set   \ \ \    } \nonumber\\
  &\boldsymbol{w}_{k,i} = \displaystyle \sum_{\ell=1}^{N}\boldsymbol{a}_{\ell
  k}(i) \boldsymbol{\psi}_{\ell,i}
\end{align}
\EndFor 
  \vspace{0.05cm}
\EndFor 
 \end{algorithmic}
 \caption{(Distributed clustering scheme)}
\label{alg:S}
\end{algorithm}

\section{Mean-Square-Error Analysis}

\noindent 
We now examine the mean-square performance of the proposed scheme.

We collect the estimates from across the network into the block vectors:
\begin{align}\label{eq:psi}
  \boldsymbol{\psi}_{i} &\triangleq\textrm{col}
\{\boldsymbol{\psi}_{1,i},\boldsymbol{\psi}_{2,i},\ldots,\boldsymbol{\psi}_{N,i}\},\\
  \boldsymbol{w}_{i} &\triangleq\textrm{col}
\{\boldsymbol{w}_{1,i},\boldsymbol{w}_{2,i},\ldots,\boldsymbol{w}_{N,i}\},\label{eq:w}
\end{align}
and define the matrices 
\begin{equation}\label{eq:AM}
{\boldsymbol{{\cal A}}}_{i}
\triangleq
{\boldsymbol{{A}}}_{i}\otimes I_M, \ \ \ {\cal M}
\triangleq\textrm{diag}
\{\mu_1,\ldots,\mu_N\}\otimes I_M.
\end{equation}
where $\boldsymbol{A}_i=[\boldsymbol{a}_{\ell k}(i)]$.
From~(\ref{eq:data_step1}) we find that the network vector $\boldsymbol{\psi}_i$ evolves over time according to 
\begin{equation}\label{eq:psi_rec_1}
   \boldsymbol{\psi}_{i}=  
   \boldsymbol{\psi}_{i-1}-{\cal M}\nabla J(\boldsymbol{\psi}_{i-1})
   -{\cal M} \boldsymbol{s}_i(\boldsymbol{\psi}_{i-1})
\end{equation}
where $\nabla J(\cdot)$ and $\boldsymbol{s}_i(\cdot)$ are
defined in~(\ref{eq:gradient}) and~(\ref{eq:gradient_noise_vec}). 
Likewise, from~(\ref{eq:data_step2}) we find that
\begin{equation}\label{eq:w_rec_1}
   \boldsymbol{w}_{i}={\boldsymbol{{\cal A}}}_{i}^\intercal
   \boldsymbol{\psi}_{i}.
\end{equation}
To proceed, we introduce the error vectors
\begin{equation}\label{eq:psi_w}
\widetilde{\boldsymbol{\psi}}_{k,i}
\triangleq
w_k^\circ-{\boldsymbol{\psi}}_{k,i}, \ \ \ 
\widetilde{\boldsymbol{w}}_{k,i}
\triangleq
w_k^\circ-{\boldsymbol{w}}_{k,i},
\end{equation}
and collect them from across the network into 
\begin{align}
\widetilde{\boldsymbol{\psi}}_i &\triangleq\textrm{col}
\{\widetilde{\boldsymbol{\psi}}_{1,i},\ldots,\widetilde{\boldsymbol{\psi}}_{N,i}\},\\
\widetilde{\boldsymbol{w}}_i
&\triangleq\textrm{col}
\{\widetilde{\boldsymbol{w}}_{1,i},\ldots,\widetilde{\boldsymbol{w}}_{N,i}\}.
\end{align}
We further define the network mean-square deviation (MSD) before and after
the fusion step  at the time  $i$ by 
 \begin{align}\label{eq:psi_rec_5}
\textrm{MSD}_{{\psi}}(i) & \triangleq 
\mathbb{E} \; \|{\widetilde{\boldsymbol{\psi}}}_{i}\|^2,\\ 
\textrm{MSD}_{{w}}(i) & \triangleq 
\mathbb{E} \; \|{\widetilde{\boldsymbol{w}}}_{i}\|^2.\label{eq:w_rec_7}
 \end{align}
\subsection{Error Dynamics}
 \noindent  Appealing to the mean-value theorem~\cite[p.~327]{book2} we can
 write
\begin{equation}\label{eq:gradient3}
\nabla J(\boldsymbol{\psi}_{i-1})=-
\boldsymbol{{\cal H}}_{i-1}\widetilde{\boldsymbol{\psi}}_{i-1}
\end{equation}
where
\begin{equation}
\boldsymbol{{\cal H}}_{i-1}
 \triangleq  \textrm{diag} \{\boldsymbol{H}_{k,i-1}\}_{k=1}^N
\end{equation}
and each matrix $\boldsymbol{H}_{k,i-1}$ is given by
\begin{equation}
\boldsymbol{H}_{k,i-1}\triangleq
\int_0^1 \nabla^2 J_k(w_k^\circ-t \widetilde{\boldsymbol{\psi}}_{k,i-1})dt.
\end{equation}
Substituting~(\ref{eq:gradient3}) into~(\ref{eq:psi_rec_1}) yields
\begin{equation}
   \boldsymbol{\psi}_{i}=  
   \boldsymbol{\psi}_{i-1}+{\cal
   M}\boldsymbol{{\cal H}}_{i-1}\widetilde{\boldsymbol{\psi}}_{i-1} -{\cal M}
   \boldsymbol{s}_i(\boldsymbol{\psi}_{i-1}).
\end{equation}
By subtracting $w^\circ$ defined in~(\ref{eq:unknown_models}) from both sides, we get
\begin{equation} \label{eq:psi_rec_3}
   \widetilde{\boldsymbol{\psi}}_{i}
=    (I_{NM}-{\cal   M}\boldsymbol{{\cal
   H}}_{i-1})\widetilde{\boldsymbol{\psi}}_{i-1}
    +{\cal M}   \boldsymbol{s}_i(\boldsymbol{\psi}_{i-1})
\end{equation}
which means that the error recursion for each individual 
agent $k$ is given by
\begin{equation}
   \widetilde{\boldsymbol{\psi}}_{k,i}
   \label{eq:psi_rec_5} =    (I_{M}-\mu_k\boldsymbol{H}_{k,i-1})
   \widetilde{\boldsymbol{\psi}}_{k,i-1}
    +\mu_k \boldsymbol{s}_{k,i}(\boldsymbol{\psi}_{k,i-1}).
\end{equation}
It is argued in~\cite[p.~347]{book2} that for step-sizes $\mu_k$ satisfying 
\begin{equation}\label{eq:mu}
0< \mu_k < \frac{2\tau_k}{ \zeta_k^2 + \beta_k^2}
\end{equation}
the mean-square-error quantity $\mathbb{E} \;
\|\widetilde{\boldsymbol{\psi}}_{k,i}\|^2$ converges exponentially  according to
the recursion:
\begin{equation}\label{eq:con_psi}
\mathbb{E} \; \|\widetilde{\boldsymbol{\psi}}_{k,i}\|^2 \leq \xi_k \ 
\mathbb{E} \; \|\widetilde{\boldsymbol{\psi}}_{k,i-1}\|^2+ \mu_k^2\rho_k^2
\end{equation}
where $0\leq \xi_k <1$ and is given by
\begin{equation}
\xi_k=1-2 \mu_k \tau_{k} +\mu_k^2(\zeta_{k}^2+\beta_k^2).
\end{equation}
It is further shown in~\cite[pp.~352,~378]{book2}  that for small step-sizes 
satisfying~(\ref{eq:mu}), the
 error recursion~(\ref{eq:psi_rec_3})
 has bounded first, second, and fourth-order moments in the following sense:
\begin{align}\label{eq:st1}
&\limsup_{i\rightarrow\infty}\; \|\mathbb{E} \; 
\widetilde{\boldsymbol{\psi}}_{i}\|
={\cal O}(\mu_{\max}) \\
&\limsup_{i\rightarrow\infty}\; \mathbb{E} \; \|
\widetilde{\boldsymbol{\psi}}_{i}\|^2
={\cal O}(\mu_{\max})\label{eq:st2} \\
&\limsup_{i\rightarrow\infty} \; \mathbb{E} \; \|
 \widetilde{\boldsymbol{\psi}}_{i}\|^4
 ={\cal O}(\mu_{\max}^2) \label{eq:st4}
\end{align}
where $\mu_{\max}$ is the maximum step-size across all agents. 

We further introduce the constant block diagonal matrix:
\begin{equation}\label{eq:h_constant}
{\cal H}\triangleq\textrm{diag}\{H_1,\ldots,H_N\}, \ \ \
{H_k}\triangleq\nabla^2J_k(w_k^\circ),
\end{equation}
and replace~(\ref{eq:psi_rec_3}) by the approximate recursion
\begin{equation}\label{eq:psi_rec_4}
   \widetilde{\boldsymbol{\psi}}'_{i}=   (I_{NM}-{\cal  
   M}{{\cal H}})\widetilde{\boldsymbol{\psi}}'_{i-1}
    +{\cal M}   \boldsymbol{s}_i(\boldsymbol{\psi}_{i-1})
\end{equation}
where the random matrix $\boldsymbol{\cal H}_{i-1}$ is replaced by ${\cal H}$. 
It was also shown in~\cite[pp.~382,~384]{book2} that, for sufficiently small
step-sizes, the error iterates that are generated by this recursion satisfy:
\begin{align}
&\lim_{i\rightarrow\infty}  \ \mathbb{E} \; 
\widetilde{\boldsymbol{\psi}}'_{i}
=0 \\
&\limsup_{i\rightarrow\infty}  \; \mathbb{E} \; \|
\widetilde{\boldsymbol{\psi}}'_{i}
\|^2 ={\cal O}(\mu_{\max}) \\
&\limsup_{i\rightarrow\infty} \; \mathbb{E} \; \|
\widetilde{\boldsymbol{\psi}}_{i}-\widetilde{\boldsymbol{\psi}}'_{i}
\|^2 ={\cal O}(\mu_{\max}^2)\\
&\limsup_{i\rightarrow\infty}\; \mathbb{E} \; \|
\widetilde{\boldsymbol{\psi}}'_{i}\|^2 =
\limsup_{i\rightarrow\infty} \; \mathbb{E} \; \|
\widetilde{\boldsymbol{\psi}}_{i}\|^2+{\cal O}(\mu_{\max}^{3/2}).
\end{align}
These results imply that, for large enough $i$, the errors
$\widetilde{\boldsymbol{\psi}}$ and $\widetilde{\boldsymbol{\psi}}'$ are
 close to each other in the mean-square-error sense. 
\subsection{One Useful Property}
  \noindent The above construction guarantees one useful property if  
  the clustering process does not incur errors of type II, 
  meaning that links that should be disconnected are indeed disconnected.
 This implies that $w_{\ell}^\circ=w_{k}^\circ$ whenever 
 $\boldsymbol{a}_{\ell k}(i)>0$.
 Using~(\ref{eq:matrix_a}), it follows that
\begin{equation}\label{eq:a_w_o_exp}
\sum_{\ell =1}^{N} \boldsymbol{a}_{\ell k}(i) w^\circ_\ell=w^\circ_k
\end{equation}
or, equivalently, 
\begin{equation}\label{eq:a_w_o}
{\boldsymbol{\cal A}}^\intercal_{i} w^\circ=w^\circ.
\end{equation}
Subtracting $w^\circ$ from both sides
of~(\ref{eq:w_rec_1}) yields,
\begin{equation}\label{eq:w_rec_4}
w^\circ-{\boldsymbol{w}}_{i}= w^\circ-{\boldsymbol{
\cal A}}^\intercal_{i}{{\boldsymbol{\psi}}}_{i}.
\end{equation}
Using~(\ref{eq:a_w_o}) we rewrite~(\ref{eq:w_rec_4}) as:
\begin{equation} \label{eq:SER}
\widetilde{\boldsymbol{w}}_{i}=
{\boldsymbol{\cal A}}^\intercal_{i}{\widetilde{\boldsymbol{\psi}}}_{i}.
\end{equation}
Taking the block maximum norm~\cite[p.~435]{book_chapter} of both sides and
using the sub-multiplicative  property of  norms implies that
\begin{equation}\label{eq:sub_2}
\|\widetilde{\boldsymbol{w}}_{i}\|_{b,\infty}
\leq \|{\boldsymbol{
\cal A}}^\intercal_{i}\|_{b,\infty} \ 
\|{\widetilde{\boldsymbol{\psi}}}_{i}\|_{b,\infty}=
\|{\widetilde{\boldsymbol{\psi}}}_{i}\|_{b,\infty}
\end{equation}
since $\boldsymbol{A}_{i}$ is left-stochastic and, therefore,
$\|\boldsymbol{\cal A}_{i}^\intercal\|_{b,\infty}=1$. It follows that
\begin{equation}\label{eq:less}
\mathbb{E}  \; \| {\widetilde{\boldsymbol{w}}}_{i}\|_{b,\infty}
\leq
 \mathbb{E}  \;
\|{\widetilde{\boldsymbol{\psi}}}_{i}\|_{b,\infty}.
\end{equation}
Results~(\ref{eq:sub_2}) and~(\ref{eq:less}) ensure that the size of the error in the $w$
domain is bounded by the size of the error in the $\psi$ domain if there are no errors of
 type II during clustering.

 \section{Performance Analysis}
 
\noindent We are ready to examine  the behavior of the probabilities of erroneous
decisions of types I and II for each agent $k$, namely, 
  the probabilities that a link between $k$ and one 
  of its neighbors will be either erroneously 
  disconnected (when it should be connected) 
  or erroneously connected (when it should be disconnected): 
\begin{align}\label{eq:type_1}
\textrm{Type-I: } w_\ell^\circ=w_k^\circ &\textrm{ and } \boldsymbol{a}_{\ell
k}(i)=0,\\
\textrm{Type-II: } w_\ell^\circ \neq w_k^\circ  &\textrm{ and }
\boldsymbol{a}_{\ell k}(i)\neq0
\label{eq:type_2}
\end{align}
for any $\ell \in {\cal N}_k$. After long enough $i$, these 
 probabilities are denoted respectively by:
\begin{align}\label{eq:P_1}
P_{\textrm{I}}&=\Pr\ (\boldsymbol{f}_{\ell k}(i)<\gamma \ | \
w_\ell^\circ= w_k^\circ),\\
P_{\textrm{II}} &=\Pr\ (\boldsymbol{f}_{\ell k}(i)\geq \gamma
\ |  \ w_\ell^\circ\neq w_k^\circ).
\label{eq:P_2}
\end{align}
Assessing the probabilities~(\ref{eq:P_1}) and~(\ref{eq:P_2}) 
is a challenging task and needs to be pursued under 
some simplifying conditions to facilitate the analysis.  This is
  due to the stochastic nature of the clustering and learning processes, 
and due to the coupling among the agents. 
Our purpose is to provide insights into the performance of 
these processes after sufficient learning time has elapsed. 
The analysis that follows adjusts the approach of~\cite{dm} to the current
setting. Different from~\cite{dm}
 where there were only two models and all agents were
trying to converge to one of these two models, 
we now have a multitude of clusters and agents that  are trying to 
converge to their own cluster model.
\subsection{Smoothing Process}
\noindent 
In order  to determine  bounds for
$P_{\textrm{I}}$ and $P_{\textrm{II}}$  we study  the probability
distribution of the trust variable $\boldsymbol{f}_{\ell k}(i)$.  
We have from~(\ref{eq:f}) that:
\begin{equation}\label{eq:f_1}
\boldsymbol{f}_{\ell k}(i)
=\nu^{i+1}\boldsymbol{f}_{\ell
k}(-1)+(1-\nu)\sum_{j=0}^i \nu^j \boldsymbol{b}_{\ell k}(i-j)
\end{equation}
where $\boldsymbol{b}_{\ell k}(i)$ is modelled as a  Bernoulli random variable with success probability $p$:
\begin{equation}\label{eq:b_1}
\boldsymbol{b}_{\ell k}(i)=
\begin{cases}
1, \textrm{with probability } p,\\
0, \textrm{with probability } (1-p).
\end{cases}
\end{equation}
We already know from~(\ref{eq:con_psi}) that, after sufficient time, the
iterates $\boldsymbol{\psi}_{k,i}$ converge to the true models
 $w_{k}^\circ$ in the mean-square-error sense.
 Hence, it is reasonable to assume that the value of 
 $p$ becomes largely time-invariant and corresponds 
 to the probability of the event described by  
\begin{equation}\label{eq:event}
\| \boldsymbol{\psi}_{\ell,i} - \boldsymbol{w}_{k,i-1}\|^2 \leq
\alpha, \textrm{ for large } i.
\end{equation}
We denote the probabilities of true and false assignments  
by
\begin{align}\label{eq:P_d_1}
P_{d}=\Pr\ (\boldsymbol{b}_{\ell k}(i)=1 \ &|\ w_\ell^\circ=
w_k^\circ),\\
P_{f}=\Pr\ (\boldsymbol{b}_{\ell k}(i)=1 \ &|  \ w_\ell^\circ \neq
w_k^\circ).
\label{eq:P_f_1}
\end{align}
These probabilities also satisfy:
  \begin{align}\label{eq:P_d_2}
 (1-P_{d})&=\Pr\ (\|\boldsymbol{\psi}_{\ell,i}-\boldsymbol{w}_{k,i-1}\|^2>\alpha
 \ | \ w^\circ_\ell=w^\circ_k),\\
  P_{f}&=\Pr\
  (\|\boldsymbol{\psi}_{\ell,i}-\boldsymbol{w}_{k,i-1}\|^2\leq\alpha \ | \
  w^\circ_\ell\neq w^\circ_k).\label{eq:P_f_2}
 \end{align}
  After  sufficient iterations, the influence of the initial
condition in~(\ref{eq:f_1}) can be ignored and we can approximate 
 $\boldsymbol{f}_{\ell k}(i)$ by the following
 random geometric series:
 \begin{equation}\label{eq:f_2} 
 \boldsymbol{f}_{\ell k}(i)\approx
 (1-\nu)\sum_{j=0}^i \nu^j \boldsymbol{b}_{\ell k}(i-j).
 \end{equation}
 As explained in~\cite{dm}, although it is generally not true, 
 we can simplify the analysis by assuming that, for large enough $i$,  
 the random variables  $\{\boldsymbol{b}_{\ell k}(m)\}$ in~(\ref{eq:f_2}) are
 independent and identically distributed. 
 This assumption is motivated by the fact that the models 
 observed by the different clusters are assumed to be sufficiently 
 distinct from each other by virtue of~(\ref{eq:delta_CR}).

Now, recall  that Markov's  inequality~\cite[p.~47]{book_pro2} implies
that for any nonnegative random variable $\boldsymbol{x}$ and positive scalar $u$, it
holds that:
\begin{equation}\label{eq:markov}
\Pr\ (\boldsymbol{x} \geq u)= \Pr\ (\boldsymbol{x}^2 \geq
u^2)\leq \frac{\mathbb{E} \; \boldsymbol{x}^2}{u^2}.
 \end{equation}
To apply~(\ref{eq:markov}) to the variable
$\boldsymbol{f}_{\ell k}(i)$, we need to determine its 
second-order moment. For
this purpose, we follow~\cite{dm} and introduce the change of variable:
 \begin{equation}\label{eq:b_2} 
 \boldsymbol{b}^\circ_{\ell k}(i-j)\triangleq  \frac{ \boldsymbol{b}_{\ell
 k}(i-j)-p} {\sqrt{p(1-p)}}.
 \end{equation}
 It can be verified that   the variables  \{$\boldsymbol{b}^\circ_{\ell
 k}(m)$\} are i.i.d.
 with zero mean and unit variance. As a result, we rewrite~(\ref{eq:f_2}) for large $i$ as:
 \begin{equation} \label{eq:f_3} 
 \boldsymbol{f}_{\ell k}(i) \approx
  p+\sqrt{p(1-p)} \boldsymbol{f}^\circ_{\ell k}(i)
 \end{equation} 
where
\begin{equation} \label{eq:f_4} 
\boldsymbol{f}^\circ_{\ell k}(i)\triangleq  (1-\nu)\sum_{j=0}^i \nu^j
\boldsymbol{b}^\circ_{\ell k}(i-j)
 \end{equation}
has  zero mean and  its variance is given by
\begin{align} 
\textrm{Var} \ [{\boldsymbol{f}^\circ_{\ell k}}(i)]
&= \mathbb{E} \; \big[{\boldsymbol{f}^\circ_{\ell k}}(i)\big]^2
-\Big[ \mathbb{E} \; {\boldsymbol{f}^\circ_{\ell k}}(i)\Big]^2 \nonumber \\ 
 &= \frac{1-\nu}{1+\nu}
(1-\nu^{2(i+1)})   \approx \frac{1-\nu}{1+\nu}.
 \end{align}
Returning to~(\ref{eq:P_1}) we now have, with $p$ replaced by $P_{d}$:
\begin{align}
P_{\textrm{I}} &\approx\textrm{Pr}
 \ (\boldsymbol{f}_{\ell k}(i)<\gamma \
| \ w_\ell^\circ= w_k^\circ)  \nonumber \\
& \leq \Pr\ \bigg( |\boldsymbol{f}^\circ_{\ell k}(i)|>
\frac{P_{d}-\gamma}{\sqrt{P_{d}(1-P_{d})}} 
\ \bigg| \ w_\ell^\circ= w_k^\circ \bigg)  \nonumber\\
&\leq \frac{1-\nu}{1+\nu}\cdot \frac{P_{d}(1-P_{d})}
{(P_{d}-\gamma)^2} \label{eq:p1}
\end{align}
where we applied~(\ref{eq:markov}) and the fact that, for any two
generic events $B_1$ and $B_2$, if $B_1$ implies $B_2$, then the probability of
event $B_1$ is less than the probability of
event $B_2$~\cite{book_pro}.  Similarly, by replacing $p$ by 
$P_{f}$, we obtain
\begin{equation}\label{eq:p2}
P_{\textrm{II}}\leq \frac{1-\nu}{1+\nu}\cdot
\frac{P_{f}(1-P_{f})}{(\gamma-P_{f})^2}.
\end{equation}
In expressions~(\ref{eq:p1}) and~(\ref{eq:p2}), 
it is assumed that the size of the threshold value
$\gamma$ used in~(\ref{eq:e}) satisfies $\gamma < P_{d}$ and $\gamma > P_{f}$. 
Since we usually desire the probability of false alarm to be small 
and the probability of detection to be close to one, 
these conditions can be met by  $\gamma\in(0,1)$. 
We show in the next section that this is indeed the case.

Results~(\ref{eq:p1}) and~(\ref{eq:p2}) 
provide bounds on the probabilities of
errors I and II. 
We next establish that   $P_{d}\rightarrow 1$ and
 $P_{f}\rightarrow 0$ to conclude
that $P_{\textrm{I}}\rightarrow 0$ and 
$P_{\textrm{II}}\rightarrow 0$.

\subsection{The Distribution of the Variables}
\noindent 
After sufficient iterations and for small enough step-sizes,
it is known that each $\boldsymbol{\psi}_{\ell,i}$ 
exhibits a distribution that is
nearly Gaussian~\cite{Opt1,Opt2,Opt3,Sacks,Nevelson,Bitmead}:
    \begin{equation}
\boldsymbol{\psi}_{\ell,i}\sim
\mathbb{N} \; (w_\ell^\circ,\mu_{\ell}  \Gamma_\ell)
 \end{equation}
where the matrix   $\Gamma_\ell$ is symmetric, positive semi-definite,
and  the solution to the following Lyapunov equation~\cite{Opt1}:
    \begin{equation}
 H_\ell
\Gamma_\ell +\Gamma_\ell
H_\ell= 
 R_{\ell}
  \end{equation}
  where the Hessian matrix $H_{\ell}$ is defined by~(\ref{eq:h_constant}) and 
   $R_{\ell}$ is the steady-state covariance matrix of the gradient noise 
  at agent $\ell$ defined by~(\ref{eq:gradient_noise_4}).
   We next introduce  the  vector
    \begin{equation}\label{eq:w_breve}
  \boldsymbol{\bar{w}}_{k,i}^\circ \triangleq
  \sum_{\ell=1}^N{\boldsymbol{a}}_{\ell k }(i)w^\circ_{\ell}.
   \end{equation}
   which should be compared with expression~(\ref{eq:data_step2}). The
   vector~(\ref{eq:w_breve}) is the result of fusing the actual models using 
   the same combination weights available at time $i$. It follows  that
   $\boldsymbol{w}_{k,i}$  exhibits a distribution 
   that is nearly Gaussian 
   since  the iterates 
   $\{\boldsymbol{\psi}_{\ell,i}\}$ can be assumed to be
    independent of each other due to the decoupled nature of their updates:
\begin{equation} 
\boldsymbol{w}_{k,i}\sim
\mathbb{N} \; ( \boldsymbol{\bar{w}}^\circ_{k,i}
,   \boldsymbol{\Omega}_{k,i})
 \end{equation}
where $\boldsymbol{\Omega}_{k,i}$ is symmetric,  positive
 semi-definite, and given by
 \begin{equation}\label{eq:Omega}
\boldsymbol{\Omega}_{k,i}\triangleq 
 \sum_{\ell=1}^N \mu_\ell \; {\boldsymbol{a}}^2_{\ell k
 }(i){\Gamma}_{\ell}.
 \end{equation}
Let 
 \begin{align}
 \boldsymbol{g}_{\ell k,i} \triangleq 
\boldsymbol{\psi}_{\ell,i}-\boldsymbol{w}_{k,i-1} 
 \end{align}  
and note  that $\boldsymbol{g}_{\ell k,i}$  is again 
   approximately Gaussian distributed
  with
\begin{equation}\label{eq:GD}
\boldsymbol{g}_{\ell k,i}\sim
\mathbb{N} \ (\boldsymbol{\bar{g}}_i, \boldsymbol{\Delta}_{\ell
k,i})
 \end{equation}
where
 \begin{equation}
\boldsymbol{\bar{g}}_i\triangleq
w^\circ_{\ell}-\boldsymbol{\bar{w}}^\circ_{k,i-1}
 \end{equation}  
 and $ \boldsymbol{\Delta}_{\ell k,i}$ is symmetric,  positive
 semi-definite, and bounded by (in view of Jensen's
 inequality\cite[p.~769]{book2}\footnote{ Since 
$\mathbb{E} \;  (\boldsymbol{a}+\boldsymbol{b})^2 \leq 2  
 \mathbb{E} \; \boldsymbol{a}^2 +
   2   \mathbb{E} \; \boldsymbol{b}^2
$.}):
   \begin{equation}\label{eq:delta}
   \boldsymbol{\Delta}_{\ell k,i} \leq 2 \; \big( \;
	\mu_\ell{\Gamma}_{\ell}+
\boldsymbol{\Omega}_{k,i-1}\;\big).
   \end{equation}
  From~(\ref{eq:Omega}),~(\ref{eq:delta}) and for any   $\ell$ and $k$ it
  holds that:
      \begin{equation}\label{eq:order}
   \boldsymbol{\Delta}_{\ell k,i} ={\cal  O}(\mu_{\max}).
   \end{equation}
    \subsection{The Statistics of
$\|\boldsymbol{g}_{\ell k,i}\|^2$}
 \noindent We now examine  the statistics of the main test variable for 
 our algorithm from~(\ref{eq:b}), namely, 
 $\|\boldsymbol{g}_{\ell k,i}\|^2$. 
 Let  $\{\mathbb{A}_{k,i-1};i>0\}$ denote the filtration 
that collects all  $\{\boldsymbol{a}_{\ell k}(i-1)\}$ 
information up to time $i-1$. Then, note that 
 \begin{align}\label{eq:mean}
\mathbb{E} \;  \left[\|\boldsymbol{g}_{\ell k,i}\|^2 
| \; \mathbb{A}_{k,i-1} \right]
&=\mathbb{E} \; \left[\textrm{Tr}\;\big(\boldsymbol{g}_{\ell k,i}
\boldsymbol{g}_{\ell k,i}^\intercal\big)\big|  \; \mathbb{A}_{k,i-1} \right] 
\nonumber\\
&=\|\boldsymbol{\bar{g}}_i\|^2+
 \textrm{Tr}  \; ( \boldsymbol{\Delta}_{\ell k,i}).
 \end{align}
Since $\boldsymbol{g}_{\ell k,i}$ is Gaussian, it holds that
  \begin{align}\label{eq:term1}
  \mathbb{E} \; & \left[ 
\|\boldsymbol{g}_{\ell k,i}\|^4 | \; \mathbb{A}_{k,i-1}  \right] \nonumber \\ 
&= 
\mathbb{E} \;  \left[  \|\boldsymbol{g}_{\ell k,i}-\boldsymbol{\bar{g}}_i+
\boldsymbol{\bar{g}}_i\|^4  | \; \mathbb{A}_{k,i-1}  \right]
\nonumber \\
&=
\mathbb{E}\; \left[  \big(
\  \|\boldsymbol{g}_{\ell k,i}-\boldsymbol{\bar{g}}_i\|^2 
+2 (\boldsymbol{g}_{\ell k,i}-\boldsymbol{\bar{g}}_i)^\intercal 
\boldsymbol{\bar{g}}_i +    \|\boldsymbol{\bar{g}}_i\|^2 
\big)^2  \big| \; \mathbb{A}_{k,i-1}  \right]  \nonumber \\
&=
\mathbb{E} \;  \left[ 
 \|\boldsymbol{g}_{\ell k,i}-\boldsymbol{\bar{g}}_i\|^4 |
 \; \mathbb{A}_{k,i-1}
\right]  \nonumber \\
&+2 \mathbb{E}   \left[  
\; \|\boldsymbol{g}_{\ell k,i}-\boldsymbol{\bar{g}}_i\|^2  
  \|\boldsymbol{\bar{g}}_i\|^2  |\; \mathbb{A}_{k,i-1}
\right]
+ \|\boldsymbol{\bar{g}}_i\|^4     
 \nonumber \\
 &+
4 \boldsymbol{\bar{g}}_i^\intercal 
\mathbb{E}\; \big[
\   (\boldsymbol{g}_{\ell k,i}-\boldsymbol{\bar{g}}_i) (\boldsymbol{g}_{\ell
k,i}-\boldsymbol{\bar{g}}_i)^\intercal | \; \mathbb{A}_{k,i-1} 
\big]\boldsymbol{\bar{g}}_i   \nonumber
\\
&= \mathbb{E} 
\; \left[  \|\boldsymbol{g}_{\ell k,i}-\boldsymbol{\bar{g}}_i\|^4  
| \; \mathbb{A}_{k,i-1}  \right]
+ 2 \; \textrm{Tr} \; ( \boldsymbol{\Delta}_{\ell k,i})
\|\boldsymbol{\bar{g}}_i\|^2 \nonumber \\
& \ + \|\boldsymbol{\bar{g}}_i\|^4
+4  \|\boldsymbol{\bar{g}}_i\|^2_{\boldsymbol{\Delta}_{\ell k,i}}
 \end{align}
where all odd order moments of $(\boldsymbol{g}_{\ell k,i}
-\boldsymbol{\bar{g}}_i)$ are zero.
 Likewise, 
  \begin{align}\label{eq:term2}
\big(\mathbb{E} \;  & \left[  \|\boldsymbol{g}_{\ell k,i}\|^2
| \; \mathbb{A}_{k,i-1}
\right]
\big)^2 \nonumber \\
&=\big(\mathbb{E} \; \left[  \|\boldsymbol{g}_{\ell
k,i}-\boldsymbol{\bar{g}}_i+\boldsymbol{\bar{g}}_i\|^2
| \; \mathbb{A}_{k,i-1}
\right]
\big)^2  \nonumber \\
&= \big( \mathbb{E}\;  \left[  \|\boldsymbol{g}_{\ell
k,i}-\boldsymbol{\bar{g}}_i\|^2 | \; \mathbb{A}_{k,i-1}
\right]+\|\boldsymbol{\bar{g}}_i\|^2
\big)^2  \nonumber \\
&=  [\; \textrm{Tr} \; ( \boldsymbol{\Delta}_{\ell k,i})\;]^2
+ 2 \; \textrm{Tr} \; ( \boldsymbol{\Delta}_{\ell
k,i})\|\boldsymbol{\bar{g}}_i\|^2 \nonumber \\
& \ +
\|\boldsymbol{\bar{g}}_i\|^4.
 \end{align}
 According  to Lemma A.2 of~\cite[p.~11]{book1},
 we have
     \begin{align}\label{eq:lemma}
\mathbb{E}  \; & \left[ 
  \|\boldsymbol{g}_{\ell k,i}-\boldsymbol{\bar{g}}_i\|^4
| \; \mathbb{A}_{k,i-1}
\right]  \nonumber \\
&= 
 [\;\textrm{Tr} \; ( \boldsymbol{\Delta}_{\ell k,i})\;]^2+
2 \; \textrm{Tr} \; ( \boldsymbol{\Delta}_{\ell k,i}^2).
  \end{align}
 Using~(\ref{eq:term1}) and~(\ref{eq:lemma}), the variance of
 $\|\boldsymbol{g}_{\ell k,i}\|^2$ is given by
 \begin{align}
\textrm{Var} \; & \left[ \|\boldsymbol{g}_{\ell k,i}\|^2 | \;
 \mathbb{A}_{k,i-1}  \right]
\nonumber \\ 
&=4 
\|\boldsymbol{\bar{g}}_i\|_{ \boldsymbol{\Delta}_{\ell k,i}}^2
+2 \; 
\textrm{Tr} \; ( \boldsymbol{\Delta}_{\ell k,i}^2).  
 \end{align}
Note from~(\ref{eq:mean}) that the  mean of $\|\boldsymbol{g}_{\ell k,i}\|^2$ is
 dominated by $\|\boldsymbol{\bar{g}}_i\|^2$ 
 for sufficiently small step-sizes. It follows from  the  Chebyshev's
 inequality~\cite[p.~455]{book_pro5} that:
    \begin{align}
 \Pr\; & \bigg( \big| \|\boldsymbol{g}_{\ell k,i}\|^2
 -\mathbb{E} \; \left[ \|\boldsymbol{g}_{\ell k,i}\|^2| \; 
 \mathbb{A}_{k,i-1} \right] \big| \geq u \bigg| \;
 \mathbb{A}_{k,i-1} \bigg) \nonumber \\
 & \leq\frac{\textrm{Var} \;  \left[ \|\boldsymbol{g}_{\ell k,i}\|^2 | \;
 \mathbb{A}_{k,i-1} \right]}{u^2}={\cal
 O}(\mu_{\max})
  \end{align}
for any constant $u>0$, which implies  that the 
 variance of $\|\boldsymbol{g}_{\ell k,i}\|^2$ is in the order of
  $\mu_{\max}$.  
    Therefore, when $w_\ell^\circ=w_k^\circ$ the probability mass of
  $\|\boldsymbol{g}_{\ell k,i}\|^2$ will concentrate around
  $\mathbb{E}(\|\boldsymbol{g}_{\ell k,i}\|^2)$, which is  in the
  order of ${\cal O}(\mu_{\max})\approx 0$. On the other hand, when 
  $w_\ell^\circ\neq w_k^\circ$, the probability mass of $\|\boldsymbol{g}_{\ell
  k,i}\|^2$ will concentrate around $\mathbb{E}(\|\boldsymbol{g}_{\ell k,i}\|^2)\approx
  \|\boldsymbol{\bar{g}}_i\|^2>0$.
   Obviously the threshold should be
 chosen as: $0<\alpha<\delta^2$, where $\delta$ is the
 clustering resolution.

 \subsection{Error Probabilities}
  
\noindent It is seen from~(\ref{eq:P_d_2}) and~(\ref{eq:P_f_2})  that
$1-P_{d}$ corresponds to the right tail probability 
of $ \|\boldsymbol{g}_{\ell k,i}\|^2$ when $w^\circ_\ell= w^\circ_k$,
and $P_{f}$ corresponds to the left 
tail probability
of  $ \|\boldsymbol{g}_{\ell k,i}\|^2$ when $w^\circ_\ell \neq w^\circ_k$.
 To examine these probabilities, we follow arguments similar
  to~\cite{cluster3} and
 apply them to the current context. We introduce the eigen-decomposition
\begin{equation}\label{eq:eigen}
\boldsymbol{\Delta}_{\ell k,i}= \boldsymbol{U}_{i}
\boldsymbol{\Lambda}_{i}
\boldsymbol{U}_{i}^\intercal
\end{equation}
where $\boldsymbol{U}_{i}$ is orthonormal and $\boldsymbol{\Lambda}_{i}$ is
diagonal and nonnegative-definite. We further introduce the normalized
variables:
  \begin{align}\label{eq:x}
\boldsymbol{x}_i
 &  \triangleq \boldsymbol{\Lambda}_{i}^{-1/2}\boldsymbol{U}_{i}^\intercal
 \boldsymbol{g}_{\ell k,i},\\
\boldsymbol{\bar{x}}_i
 & \triangleq \boldsymbol{\Lambda}_{i}^{-1/2}\boldsymbol{U}_{i}^\intercal
 \boldsymbol{\bar{g}}_{i}
 \label{eq:xbar}
 \end{align}
and it follows from~(\ref{eq:GD}),~(\ref{eq:x}), and~(\ref{eq:xbar}) that 
  \begin{equation}\label{eq:GDX}
\boldsymbol{x}_i\sim
\mathbb{N} \ (\boldsymbol{\bar{x}}_i,
\mathds{I}_M).
\end{equation}
Note also from~(\ref{eq:x}) that
\begin{equation}\label{eq:11}
\|\boldsymbol{g}_{\ell k,i}\|^2=\boldsymbol{x}_i^\intercal \boldsymbol{\Lambda}_{i}
\boldsymbol{x}_i =\sum_{h=1}^M \boldsymbol{\lambda}_{h,i}
\boldsymbol{x}_{h,i}^2.
\end{equation}
where  $\boldsymbol{x}_{h,i}$ denotes the $h-$th element of 
$\boldsymbol{x}_i$
and $\boldsymbol{\lambda}_{h,i}$ denotes the $h-$th
 diagonal element of $
\boldsymbol{\Lambda}_{i}$.
 \subsubsection{The probability $1-P_{d}$}
\noindent It follows from the inequality 
\begin{equation}\label{eq:12}
\|\boldsymbol{g}_{\ell k,i}\|^2=\boldsymbol{x}_i^\intercal \boldsymbol{\Lambda}_{i}
 \boldsymbol{x}_i
\leq \|\boldsymbol{\Delta}_{\ell k ,i}\| \cdot \|\boldsymbol{x}_i\|^2,
\end{equation}
that the following relation is satisfied
\begin{equation}\label{eq:13}
\{\|\boldsymbol{g}_{\ell k,i}\|^2> \alpha\} \subseteq
\{ \|\boldsymbol{\Delta}_{\ell k ,i}\| \cdot \|\boldsymbol{x}_i\|^2 >
\alpha\}.
\end{equation}
Defining 
\begin{equation}\label{eq:alpha_bold}
\boldsymbol{\alpha}_k(i)\triangleq \alpha/ \|\boldsymbol{\Delta}_{\ell k ,i}\|
\end{equation},
we can write
using~(\ref{eq:13}):
  \begin{equation}\label{eq:P_d_3}
\Pr\ (\|\boldsymbol{g}_{\ell k,i}\|^2>\alpha
 \ | \ w^\circ_\ell=w^\circ_k) \leq 
 \Pr\ ( \|\boldsymbol{x}_i\|^2>\boldsymbol{\alpha}_k(i)
 \ | \ \boldsymbol{\bar{x}}_i=0).
 \end{equation}
We know from~(\ref{eq:GDX}) that 
\begin{equation}\label{eq:ChiD}
\|\boldsymbol{x}_i\|^2\sim {\cal X}_M^2
\end{equation}
where $ {\cal X}_M^2$ denotes the Chi-square distribution with $M$
degrees of freedom and its  mean value is $M$.
According to the Chernoff bound for the central Chi-square 
distribution with $M$ degrees of
freedom\footnote{Let 
$\boldsymbol{y}\sim {\cal X}_r^2$. Acoording to the Chernoff bound for the
central Chi-square distribution with $r$ degrees of
freedom, for any  $\epsilon >0$ it holds
that~\cite[p.~2501]{Chernoff}: $ \Pr\; 
(\boldsymbol{y}>r(1+\epsilon))\leq   \textrm{exp} \; [-\frac{r}{2} (\epsilon - \textrm{log}
(1+\epsilon))]$.} we have
  \begin{align}\label{eq:PB1}
 \Pr\; & \big( \; \|\boldsymbol{x}_i\|^2>\boldsymbol{\alpha}_k(i)
 \; \big| \; \boldsymbol{\bar{x}}_i=0 \; \big) \nonumber \\
 &=  \Pr\;  \bigg( \|\boldsymbol{x}_i\|^2>
 \frac{M \cdot \boldsymbol{\alpha}_k(i)}{M}
 \ \bigg| \ \boldsymbol{\bar{x}}_i=0 \bigg)  \nonumber \\
 &\leq \textrm{exp} \Big[- \frac{M}{2}\Big( 
  \frac{\boldsymbol{\alpha}_k(i)}{M}
  -\textrm{log}\;\Big(1+\frac{\boldsymbol{\alpha}_k(i)}{M}-1\Big)
  \Big)  \Big]  \nonumber \\
  &= \Big(\frac{\boldsymbol{\alpha}_k(i) \cdot e}{M}\Big)^{M/2} \cdot
  \textrm{exp} \Big[-\frac{\boldsymbol{\alpha}_k(i)}{2}  \ \Big]
 \end{align}
 where $e$ is
 Euler's number.  
 For
 small enough step-sizes we conclude 
 from~(\ref{eq:order}),~(\ref{eq:alpha_bold}), and~(\ref{eq:PB1}) that after
 sufficient iterations, it holds that:
 \begin{align}\label{eq:pd}
 (1-{P}_{d})
\leq {\cal O}(e^{-c_1/\mu_{\max}})
 \end{align}
 for some constant $c_1>0$. 
 
\subsubsection{The probability $P_{f}$}
\noindent The approximate characteristic function of 
$\|\boldsymbol{g}_{\ell k,i}\|^2$~\cite[Eq.~(118)]{cluster3} when
$w^\circ_\ell\neq w^\circ_k$ is given by:
\begin{equation}
c_{\|\boldsymbol{g}_{\ell k,i}\|^2}(t)\approx
e^{jt \|w^\circ_\ell-w^\circ_k\|^2 - 2t^2
\|w^\circ_\ell-w^\circ_k\|^2_{\boldsymbol{\Lambda_i}}}
\end{equation}
which implies that for sufficiently small $\mu_{\max}$,
\begin{equation}\label{eq:GDr}
\|\boldsymbol{g}_{\ell k,i}\|^2\sim
\mathbb{N} \; (\|w^\circ_\ell-w^\circ_k\|^2,4
\|w^\circ_\ell-w^\circ_k\|^2_{\boldsymbol{\Lambda_i}}).
\end{equation}
Therefore, from~\cite{cluster3}\footnote{Let 
$\boldsymbol{y}\sim \mathbb{N} \; (0,1)$. 
Acoording to the Chernoff bound 
for the Gaussian error function it holds
that~\cite{Chernoff2}: \; $Q(\boldsymbol{y})\leq \frac{1}{2} \textrm{exp}
\; [-\frac{\boldsymbol{y}^2}{2}]$.} we obtain that
  \begin{align}\label{eq:PB2}
 \Pr& \ (\|\boldsymbol{g}_{\ell k,i}\|^2<\alpha
 \ | w^\circ_\ell\neq w^\circ_k \ )
 \approx  { Q}\bigg(
 \frac{ \|w^\circ_\ell- w^\circ_k\|^2- \alpha}
 {2
 \|w^\circ_\ell-w^\circ_k\|_{\boldsymbol{\Lambda_i}}}
 \bigg)  \nonumber \\
 & \leq \frac{1}{2} \;\textrm{exp}
\; \bigg[-\frac{(\|w^\circ_\ell-w^\circ_k\|^2-
 \alpha)^2}{8 
 \|w^\circ_\ell-w^\circ_k\|^2_{\boldsymbol{\Lambda_i}}}
 \bigg]
 \end{align}
where the letter $Q$ refers here  to the traditional  $Q-$function (the tail probability of the
 standard Gaussian distribution). For small 
 enough step-sizes, after sufficient iterations and
 from~(\ref{eq:order}),~(\ref{eq:eigen}), and~(\ref{eq:PB2}), it holds that
  \begin{align}\label{eq:pf}
 {P}_{f}
 &\leq {\cal O}(e^{-c_2/\mu_{\max}})
 \end{align}
 for some constant $c_2>0$. It is then seen that the
  probabilities $P_{\textrm{I}}$
 and $P_{\textrm{II}}$ are expected to
  approach zero exponentially fast for  vanishing step-sizes.

\section{Linking Application}

\subsection{Clustering With Linking Scheme}
\begin{figure}
\centering\includegraphics[width=17pc]{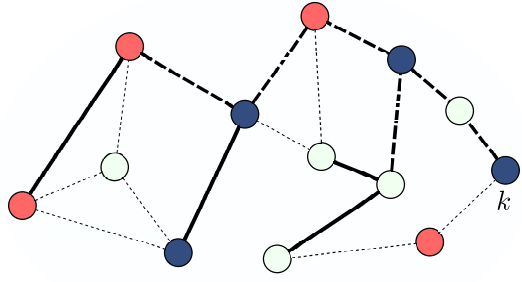}
\caption{\small The clustered and linked topology  that will result for the
network shown in Fig.~\ref{fig:Partition0}. The bold dashed  lines depict the
links used for relaying data among agents.}
\label{fig:Partition2}
\end{figure}

\noindent We propose in this section an additional mechanism to enhance the
performance of each cluster by using  the unused links  to relay 
information. Figure~\ref{fig:Partition2} shows the linked topology 
that results for the same example shown earlier in Fig.~\ref{fig:Partition0}(b). 
The figure shows that the links which are supposed to be unused for
sharing data among neighbors belonging to different clusters,
are used now to relay data among agents.

We assume in this section that the links among agents are
\emph{symmetric}, i.e. if $\ell \in {\cal N}_k \Longleftrightarrow k \in
{\cal N}_\ell$. Under normal operation, each agent $k$ will be
 receiving and processing iterates only from those neighbors that 
 it believes belong 
to the same cluster as $k$.

 We modify this operation by allowing $k$ 
to receive iterates from {\em all} of its neighbors. 
It will continue to use the iterates from neighbors in the same cluster to
update its weight estimate $\boldsymbol{w}_{k,i}$.
 The iterates that arrive from
neighbors that may belong to other clusters are not
 used during this fusion process. 
 Instead, they will be relayed forward by agent $k$ as follows. 
 For each of its neighbors $\ell\in{\cal N}_k$, 
 agent $k$ will send  $\boldsymbol{\psi}_{k,i}$
  and  another 
  vector $\boldsymbol{\phi}_{k \ell,i}$ 
 The vector $\boldsymbol{\phi}_{k\ell,i}$  is constructed as follows. Agent $k$ 
  chooses from among all the iterates it receives  from its neighbors, 
  that iterate that is closest to $\boldsymbol{\psi}_{\ell,i}$:
\begin{equation}\label{eq:phi}
\boldsymbol{\phi}_{k \ell,i} \; =\ \;
\underset{\hspace{-2.6cm} \substack{\{k,m\} \\ \forall m\in{\cal
N}_k,\;m\notin {\cal N}_{\ell}}}
{\mbox{\textrm arg min}  \ \  \|\boldsymbol{\psi}_{m,i}-
\boldsymbol{\psi}_{\ell,i}\|^2.}
\end{equation}
Observe that the minimization is over $k$ and all neighbors of $k$ that are not
neighbors of $\ell$. This condition is important  to avoid  receiving the same
information multiple times. Observe also that  under this scheme, agent $k$
will need to receive the iterates from all of its neighbors 
(those that it believes belong to its clusters 
and those that do not); it also needs to receive information 
about their neighborhoods, i.e., the ${\cal N}_{\ell}$ for 
each of its neighbors $\ell$.

 The following steps describe the clustering with linking algorithm.
 We collect all $\{\boldsymbol{\phi}_{\ell k,i}\}$ into  a matrix
 $\boldsymbol{\Phi}_i$.
By setting $\gamma=0.5$ in Eq.~(\ref{eq:e}) the operation of setting each entry
$\boldsymbol{e}_{\ell k}(i)$ becomes rounding to the nearest integer and is
denoted by $\lfloor \cdot\rceil$.

%
  \begin{algorithm}[H]
  \begin{algorithmic}
   \State Initialize $\boldsymbol{F}_{-1}=\boldsymbol{B}_{-1}=\boldsymbol{E}=I$,
   $\boldsymbol{\Phi}_{-1}=0$, and
   $\boldsymbol{\psi}_{-1}=\boldsymbol{w}_{-1}=0$.
    \For {$i\geq 0$}
      \vspace{0.05cm}
   \For {$k=1,\ldots,N$}
    \vspace{-0.12cm}
       \begin{align}
   & \ \boldsymbol{\psi}_{k,i}= \boldsymbol{\psi}_{k,i-1}- \mu_k \widehat{\nabla
       J_k}(\boldsymbol{\psi}_{k,i-1})  \ &
  \end{align}
    \For {$\ell \in {\cal N}^-_{k}$}
 \begin{align} 
 &\ \ \  \ \ \ \ \ \textrm{send }  \boldsymbol{\psi}_{k,i} \textrm{ and } 
 \boldsymbol{\phi}_{k
 \ell ,i-1} \nonumber \\ 
 &\ \ \  \ \ \  \ \ \textrm{receive }  \boldsymbol{\psi}_{\ell,i} \textrm{ and } 
 \boldsymbol{\phi}_{\ell k
 ,i-1} \nonumber \\ 
 &\ \ \  \ \ \  \ \ \boldsymbol{b}_{\ell k}(i)=\begin{cases}
 1, &   \textrm{if } \| \boldsymbol{\phi}_{\ell k,i-1} -
 \boldsymbol{w}_{k,i-1}\|^2 \leq \alpha\\
 0,  & \textrm{otherwise}
 \end{cases}\\
 &\ \ \  \  \ \ \ \ \boldsymbol{f}_{\ell k}(i)=\nu  
 \boldsymbol{f}_{\ell k}(i-1)+ (1-\nu)  \boldsymbol{b}_{\ell
 k}(i)\\
 &\ \ \  \  \  \ \ \  \boldsymbol{e}_{\ell k}(i)=\lfloor \boldsymbol{f}_{\ell
 k}(i)\rceil
\end{align}
\EndFor 
\begin{align}
    & \textrm{select } \{\boldsymbol{a}_{\ell k}(i)\} \textrm{ according
    to}~(\ref{eq:matrix_a}) \textrm{ and set \ \  } \nonumber
    \\
    &  \boldsymbol{w}_{k,i}= \sum_{\ell=1}^{N}\boldsymbol{a}_{\ell k}(i)
 \boldsymbol{\phi}_{\ell k,i-1}\\
& \textrm{update } \{\boldsymbol{\phi}_{k \ell}(i)\} \textrm{ according
    to}~(\ref{eq:phi}) \nonumber
      \end{align}
 \EndFor
  \EndFor
 \end{algorithmic}
 \caption{(Clustering with linking scheme)}
\label{alg:S}
  \end{algorithm} 

\section{Simulation Results}
\begin{figure}
 \centering 
 \psfrag{V}[b]{\small $\sigma_{v,k}^2$}
  \psfrag{R}[b]{\small$\textrm{Tr}(R_{u,k})$}
   \centerline{\includegraphics[width=20pc]{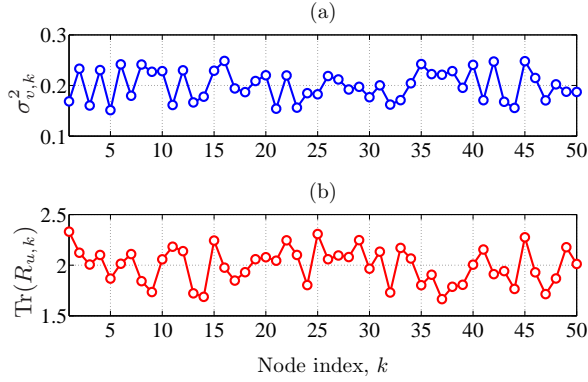}}
   \caption{ \small The statistical noise and signal profiles over the network.}
\label{fig:CNoise}
\end{figure}

 \noindent  We consider a fully connected network with 50 randomly distributed 
  agents.
The agents observe data originating from three different models ($C=3$).
   Each model $w_{{\cal C}_m}^\circ\in \mathbb{R}^{M\times1}$ is generated as
   follows: $ w^\circ_{\mathcal{C}_m}=[ w_{r_1},\ldots,w_{r_M}]^\intercal$,
   with entries     $w_{r_c}\in[1,-1]$.
    In our example we set $M=2$; larger values of $M$
     are generally easier for clustering and, therefore,    we illustrate the operation of the algorithm for $M=2$. 
   The assignment of the agents to models is random.
   Agents having the same color belong to the same cluster. The maximum
   number of neighbors is $n_{\max}= 6$. Every agent $k$ has
    access to a scalar measurement $\boldsymbol{d}_k(i)$ 
   and a $1\times M$ regression vector
 $\boldsymbol{u}_{k,i}$. The measurements across the agents are assumed
 to be generated via the  linear regression model
$
 \boldsymbol{d}_k(i)=\boldsymbol{u}_{k,i}{w}_{k}^\circ+\boldsymbol{v}_{k}(i)
$, where $\boldsymbol{v}_{k}(i)$ is measurement noise assumed to be a zero-mean
 white random process that is  independent over space. It is also assumed
 that the regression data $\boldsymbol{u}_{k,i}$ is independent over space and
 independent of $\boldsymbol{v}_{\ell}(j)$ for all $k,\ell,i,j$. All random
 processes are assumed to be stationary. The statistical profile of
  the noise across the agents for $k=1,\ldots,N$ is shown in
  Fig.~\ref{fig:CNoise}(a). The regressors are of size $M=2$ 
  and have diagonal
  covariance matrices ${R}_{u,k}$ shown in Fig.~\ref{fig:CNoise}(b).
  We set $\{\mu,\alpha,\nu,\delta,\gamma\}=\{0.05,0.015,0.98,0.17,0.5\}$.  
  We use the uniform combination policy to generate the coefficients
   $\{a_{\ell
 k}(i)\}$.
     \begin{figure}
 \vspace{-0cm}
 \centering 
   \centerline{\includegraphics[width=22pc]{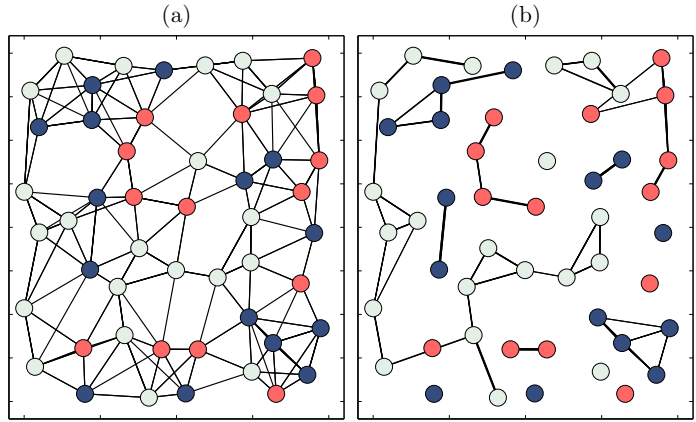}}
   \caption{ \small The network topology (a) and the clustered  topology at 
   steady-state (b).}
\label{fig:CTopology}
\end{figure}

\begin{figure}
 \vspace{-0cm}
 \psfrag{simulated-Zhao}[bl]{ \scriptsize simulated \cite{cluster3}}
  \psfrag{simulated-Chen-CTA}[bl]{ \scriptsize simulated \cite{cluster4}}
   \psfrag{MSD-w}[bl]{ \scriptsize simulated (MSD)}
 \centering 
   \centerline{\includegraphics[width=22pc]{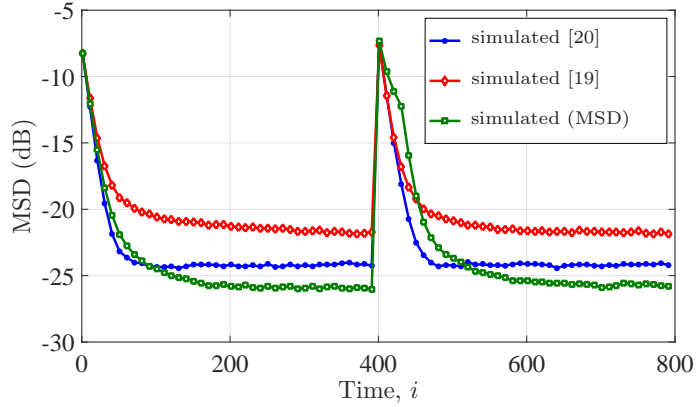}}
   \caption{ \small The transient mean-square deviation (MSD) using different
   approaches.}
\label{fig:CMSD}
\end{figure}

Figure~\ref{fig:CTopology}(a) shows the topology of one of
   100 Monte Carlo experiments.
   Figure~\ref{fig:CTopology}(b) presents the
   final topology after applying the  clustering technique.  
 Figure~\ref{fig:CMSD}(a)
 depicts the simulated transient mean-square deviation (MSD) of the network 
 compared to other
 clustering methods.  The 
 model assignments change at time instant $i=400$. 

\begin{figure}
 \vspace{-0cm}
 \psfrag{Chen-CTA}[bl]{ \scriptsize  \cite{cluster4} scheme}
  \psfrag{Zhao}[bl]{ \scriptsize  \cite{cluster3} scheme}
   \psfrag{me}[bl]{ \scriptsize proposed scheme}
 \psfrag{v1}[bl]{ \small $\overline{v}_{\textrm{I}}$}
  \psfrag{v2}[bl]{ \small $\overline{v}_{\textrm{II}}$}
 \centering 
   \centerline{\includegraphics[width=23pc]{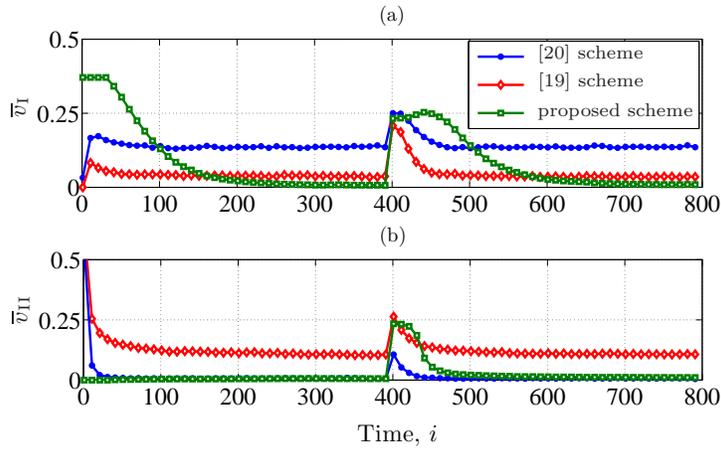}}
   \caption{ \small The normalized clustering errors of types I and II over the
   network.}
\label{fig:CErrors}
\end{figure}

  The normalized clustering
 errors of types I and II  by each agent $k$ at time  $i$ are given, 
 respectively, by
\begin{align}\label{eq:em}
\boldsymbol{v}_{\textrm{I},k}(i)
&\triangleq\frac{(\mathds{1}-{[\boldsymbol{E}_i]}_{:,k})^\intercal
\times({[E^\circ]}_{:,k}-{[\boldsymbol{E}_i]}_{:,k})}{(n_{k}-1)}\\
\boldsymbol{v}_{\textrm{II},k}(i)
&\triangleq\frac{{[\boldsymbol{E}_i]}_{:,k}^\intercal\times({[\boldsymbol{E}_i]}_{:,k}-{[E^\circ]}_{:,k})}{(n_{k}-1)}
\label{eq:ef}
\end{align}
 where  $E^\circ$ is the true clustering
 matrix. Figures~\ref{fig:CErrors}(a)--\ref{fig:CErrors}(b) depict the
 normalized clustering errors
  $\overline{v}_{\textrm{I}}$ and $\overline{v}_{\textrm{II}}$ over
   the network.


\begin{figure}
 \centering 
   \vspace{0.1cm}
    \centerline{\includegraphics[width=22pc]{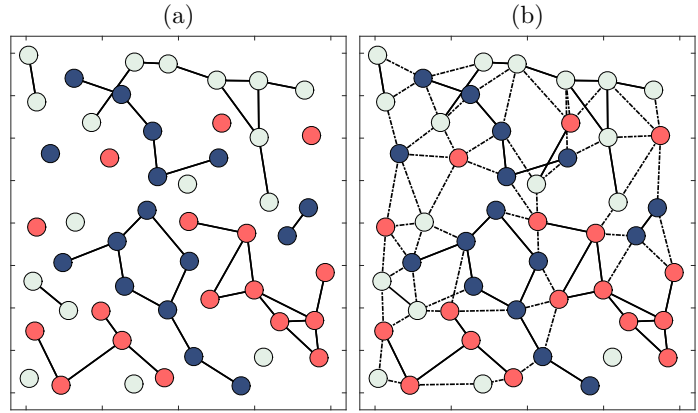}}
   \caption{ \small The clustered network topology at 
   steady-state (a) and the clustered and
   linked   topology at 
   steady-state (b).}
\label{fig:LTopology}
\end{figure}

\begin{figure}
  \vspace{0cm}
 \psfrag{v1}[bl]{ \scriptsize $v_{\textrm{I}}$}
  \psfrag{v2}[bl]{ \scriptsize $v_{\textrm{II}}$}
   \psfrag{simulated with linking}[bl]{ \hspace{-0.02cm} \scriptsize clustering
   with linking}
  \psfrag{simulated without linking}[bl]{\hspace{-0.02cm} \scriptsize clustering
  without linking}
   \centerline{\includegraphics[width=26pc]{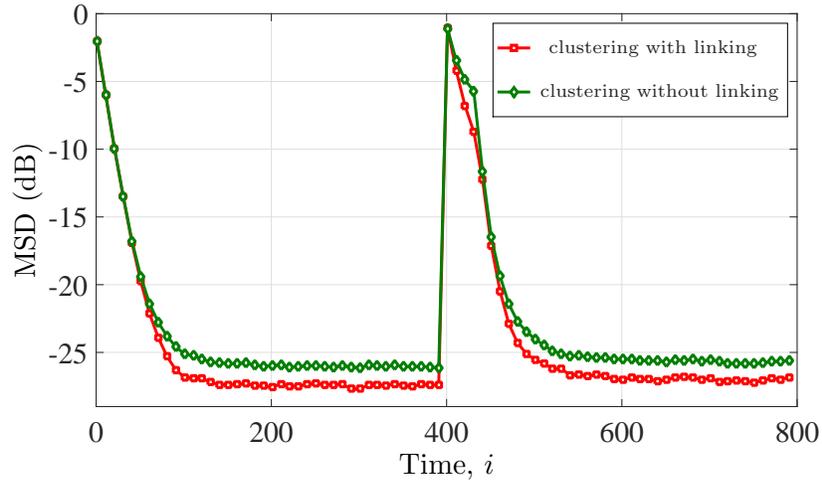}}
   \caption{ \small The transient mean-square deviation
   with and without applying the linking technique.}
\label{fig:LMSD}
\end{figure}

\begin{figure}
   \vspace{0cm}
 \psfrag{v1}[bl]{ \scriptsize $\overline{v}_{\textrm{I}}$}
  \psfrag{v2}[bl]{ \scriptsize $\overline{v}_{\textrm{II}}$}
   \psfrag{simulated with linking}[bl]{ \scriptsize Simulated with linking}
  \psfrag{simulated without linking}[bl]{\scriptsize Simulated without linking}
   \centerline{\includegraphics[width=22pc]{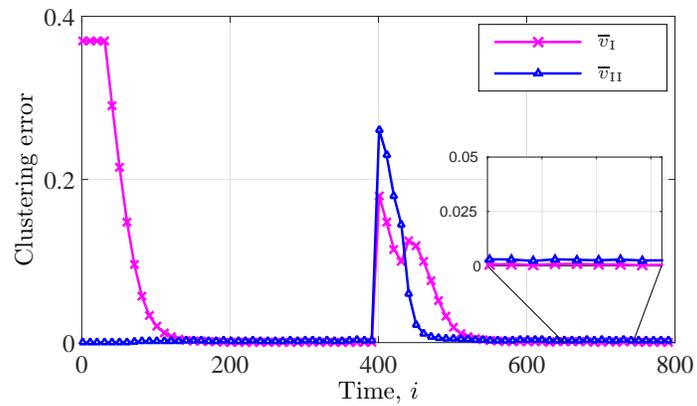}}
   \caption{ \small The normalized
   clustering errors of types I and II over the network.}
\label{fig:LErrors}
\end{figure}

   Using the same setup of the previous example,
  Fig.~\ref{fig:LTopology}(a) shows the topology of 
  one  experiment with  the clustering technique only.  
  Figure~\ref{fig:LTopology}(b) presents 
  the final topology when  we apply the
  clustering with linking technique.
  Figure~\ref{fig:LMSD} indicates the simulated  transient mean-square 
  deviation
  (MSD) of the agents with and without the linking technique.
   The normalized clustering errors over the network are shown in
 Fig.~\ref{fig:LErrors}.

\section{ Conclusion}
 \noindent We  proposed a distributed algorithm 
 that carries  out the tasks of estimation and 
  clustering simultaneously with
 exponentially decaying error probabilities for false decisions.
  We  showed how the agents
    choose the  subset of their neighbors to cooperate with and turn off 
    suspicious links.
     The simulations illustrate
    the performance of the proposed strategy and compare with other
     related works.
      We  proposed an additional step to enhance the
     performance by linking, as much as possible, 
      the agents that belonging to 
     the same cluster and do not have
     direct links to connect them.
    
  \bibliographystyle{plain}
  \bibliography{D.bib}

\end{document}